\documentclass[12pt,twoside,reqno]{amsart}

\newtheorem{thm}{Theorem}[section]
\newtheorem{lemma}{Lemma}[section]
\newtheorem{remark}{Remark}[section]
\newtheorem{proposition}{Proposition}[section]
\newtheorem{corollary}{Corollary}[section]

\def \R{{\Bbb R}}
\def \N{{\Bbb N}}

\numberwithin{equation}{section}

\begin{document}

\title[Modified Zakharov-Kuznetsov equation]
{ Modified Zakharov-Kuznetsov equation on rectangles }
\author[
M. Castelli, \ G. Doronin] { M. Castelli$^\ast$, \ G. Doronin
%$^\S$
%$^\ast$
\bigskip
\\
{\tiny
Departamento de Matem\'atica,\\
Universidade Estadual de Maring\'a,\\
87020-900, Maring\'a - PR, Brazil. } }
\address
{
Departamento de Matem\'atica\\
Universidade Estadual de Maring\'a\\
87020-900, Maring\'a - PR, Brazil. }
\email{marcos\_castelli@hotmail.com \ \ ggdoronin@uem.br}
%\date{}

\subjclass {35M20, 35Q72} \keywords {mZK equation, well-posedness}
%\thanks{$^\ast$Supported by ....}
\thanks{$^\ast$Partially supported by CAPES}
%\thanks{$^\S$Corresponding author, partially supported by Funda\c{c}\~ao Arauc\'aria, Pr., Brazil}

\begin{abstract}
Initial-boundary value problem for the modified Zakharov-Kuznetsov
equation posed on a bounded rectangle is considered. The main
difficulty is the critical power in nonlinear term. The results on
existence, uniqueness and asymptotic behavior of solutions are
presented.
\end{abstract}

\maketitle

\section{Introduction}\label{introduction}

We are concerned with initial-boundary value problems (IBVPs) posed
on bounded rectangles located at the right half-plane
$\{(x,y)\in\mathbb{R}^2:\ x>0\}$  for the modified
Zakharov-Kuznetsov (mZK) equation \cite{pastor1}
\begin{equation}
u_t+u_x+u^2u_x+u_{xxx}+u_{xyy}=0.\label{mzk}
\end{equation}
This equation is a generalization \cite{pastor} of the classical
Zakharov-Kuznetsov (ZK) equation \cite{zk} which is a
two-dimensional analog of the well-known modified Korteweg-de Vries
(mKdV) equation \cite{bona2}.

Note that both ZK and mZK possess real plasma physics applications
\cite{zk,mzk}.

As far as ZK is concerned, the results on both IVP and IBVPs can be
found in
%readers who are interested in both IVP and IBVP are referred to
\cite{faminski,faminski2,farah,larkin,pastor,pastor2,saut,rosier2,temam,
temam2}. For IVP to mZK, see \cite{pastor1}; at the same time we do
not know solid results concerning IBVP to mZK. The main difference
between initial and initial-boundary value problems is that IVP on
$\mathbb{R}^2$ provides (almost immediately) good estimates in
$H^1({\mathbb{R}^2})$ by the conservation laws \cite{pastor1}, while
IBVP does not possesses this advantage.

Our work is motivated by \cite{temam} and provides a natural
continuation of \cite{doronin} where the original ZK equation was
considered. There one can find out a more detailed background,
descriptions of main features, and the deployed reference list.

In the present note we put forward an analysis of \eqref{mzk} posed
on a bounded rectangle with homogeneous boundary conditions. Since
the power is critical \cite{pastor,pastor1}, a challenge concerning
the well-posedness of IBVPs appears. Section \ref{local existence}
provides the local results via fixed point arguments. In Section
\ref{global} we obtain global estimates which simultaneously provide
the exponential decay rates of solution. These results have been
proven for sufficiently small initial data, and under domain's size
restrictions. Restrictions upon the domain appear naturally due to
the presence of a linear transport term $u_x,$ see
\cite{doronin,rosier} for details. For one-dimensional dispersive
models the critical nonlinearity has been treated in
\cite{larkin19}.

\section{Problem, notations and preliminaries}\label{problem}

%$\overset{c}{\hookrightarrow}$

Let $L,B,T$ be finite positive numbers. Define $\Omega$ and $Q_T$ to
be spatial and time-spatial domains
\begin{equation*}
\Omega=\{(x,y)\in\mathbb{R}^2: \ x\in(0,L),\ y\in(-B,B) \},\ \ \
Q_T=\Omega \times (0,T).
\end{equation*}

In $Q_T$ we consider the following IBVP:
\begin{align}
P&u\equiv u_t+u_x+u^2u_x+u_{xxx}+u_{xyy}=0,\ \ \text{in}\ Q_T;
\label{2.1}
\\
&u(x,-B,t)=u(x,B,t)=0,\ \ x\in(0,L),\ t>0;
\label{2.2}
\\
&u(0,y,t)=u(L,y,t)=u_x(L,y,t)=0,\ \ y\in(-B,B),\ t>0;
\label{2.3}
\\
&u(x,y,0)=u_0(x,y),\ \ (x,y)\in\Omega, \label{2.4}
\end{align}
where $u_0:\Omega\to\mathbb{R}$ is a given function.

Hereafter subscripts $u_x,\ u_{xy},$ etc. denote the partial
derivatives, as well as $\partial_x$ or $\partial_{xy}^2$ when it is
convenient. Operators $\nabla$ and $\Delta$ are the gradient and
Laplacian acting over $\Omega.$ By $(\cdot,\cdot)$ and $\|\cdot\|$
we denote the inner product and the norm in $L^2(\Omega),$ and
$\|\cdot\|_{H^k}$ stands for the norm in $L^2$-based Sobolev spaces.
Abbreviations like $(L^s_t;L^l_{xy})$ are also used for anisotropic
spaces.

To prove the results we will apply
\begin{lemma}  \label{prop 1 decaimento} \text{\sl (V. A. Steklov)}
Let $L,B>0$ and $\omega \in H_0^1(\Omega)$. Then
\begin{equation}\label{decaimento 1.4}
\int_0^L \int_{-B}^B \omega^2(x,y)dxdy \leq
\frac{4B^2}{\pi^2}\int_0^L \int_{-B}^B \omega_y^2(x,y)dxdy,
\end{equation}
%-------------------------------
%-------------------------------
and
\begin{equation}\label{decaimento 1.5}
\int_0^L \int_{-B}^B \omega^2(x,y)dxdy \leq
\frac{L^2}{\pi^2}\int_0^L \int_{-B}^B \omega_x^2(x,y)dxdy.
\end{equation}
\end{lemma}
See \cite{doronin} for the proof.

The Nirenberg theorem (also often called as the Gagliardo-Nirenberg
inequality) will be used in the following form:
\begin{lemma}\label{prop 2 decaimento} \text{\sl (L. Nirenberg)}
Let $U\subset \mathbb{R}^n$ be a bounded open set, $u \in L^q(U)$
and $D^m u\in L^r(U)$, $1\leq q , r \leq \infty .$ For $0\leq j \leq
m ,$ the following inequality holds:
\begin{equation}\label{dsg Nrbg}
\|D^j u\|_{L^p(U)}  \leq C_{U} \|D^m u\|_{L^r(U)} ^{\alpha} \|
u\|_{L^q(U)} ^{1-\alpha}
\end{equation}
where
$$
\frac{1}{p}= \frac{j}{n}+\alpha\left(\frac{1}{r}-\frac{m}{n}\right) + (1-\alpha)\frac{1}{q},
$$
for all $\alpha$ from the interval $$ \frac{j}{m}\leq \alpha \leq
1.$$ The constant $C_{U}$ depends on $n,m,j,q,r,\alpha .$
\end{lemma}
For the proof see \cite{evans}.

\begin{corollary}\label{corollary Nirenberg} Let $u \in H_0^1(\Omega).$ Then for all $p \in {\N}$
$$
\| u\|_{L^{2p}(\Omega)}  \leq C_{2p} \| \nabla u\| ^{\frac{p-1}{p}} \| u\| ^{\frac{1}{p}}
$$
where
$$
C_{2p}= \left(\frac{ p!}{\sqrt{2}^{p-1}}\right)^\frac{1}{p}.
$$
\end{corollary}
The result can be proved by induction. We will also use the simple
\begin{lemma}\label{prop 2 decaimento} Let $u\in H^1(\Omega)$ and $u_{xy} \in L^2(\Omega).$
Then
$$
\sup_{(x,y) \in \Omega} u^2(x,y) \leq \|u\|_{H^1(\Omega)}^2 + \|u_{xy}\|_{L^2(\Omega)}^2.
$$
\end{lemma}
See \cite{doronin} for the proof.

\section{Existence in sub-critical case}\label{existence}

In this section we state the existence result in sub-critical case,
i.e., for $\delta\in(0,1).$ Technically, we mainly follow
\cite{temam2}. A short motivation for this study is provided in
subsection \ref{motivation}.

\section{Local results}\label{local existence}

Consider the following Cauchy problem in abstract form:
\begin{equation}\label{Eq sem. grp.}
 \left\{
\begin{array}{c}
%-------------------------------
u_t + Au= f, \;\; \\ %\text{em} \;\;  \mathrm{Q}_T = \Omega\times (0,T), \;\; ,\\
%-------------------------------
%\hspace{-0.5cm}u=0\;\;\text{em}\;\;\partial \Omega , \;\; u_x(L,y,t)=0 , \\
%-------------------------------
\hspace{0cm}u(0)=u_0 ,\;\;%\text{em}\;\;\Omega ;
%-------------------------------
\end{array}
\right.
\end{equation}
where $f\in L^1(0,T;L^2(\Omega))$ and $A: L^2(\Omega)\to
L^2(\Omega)$ defined as $A\equiv
\partial_x + \Delta\partial_x $ with the domain
$$
D(A) = \{ u \in L^2(\Omega ) \, ; \, \Delta u_x + u_x \in
L^2(\Omega) \, ;\ u|_{\partial \Omega}=0 \ \text{ and }\ u_x(L,y,t)
= 0, \ t\in (0,T) \},
$$
endowed with its natural Hilbert norm $ \|u\|_{D(A)}(t) = \left(
\|u\|^2_{L^2(\Omega)}(t) + \| \Delta u_x + u_x \|^2_{L^2(\Omega)}(t)
\right)^{1/2} $ for all $t\in (0,T)$.
%\begin{proposition}\label{proposition1} $D(A) \hookrightarrow  H^1_0(\Omega)\cap H^2(\Omega)$ com imersao continua.
%\end{proposition}
%Prova em \cite{temam2}.
\begin{proposition}\label{res. 1 teman} Assume $u_0\in D(A)$ and $f\in L^1_{loc}(\R^+;L^2(\Omega))$
with $f_t\in L^1_{loc}(\R^+;L^2(\Omega))$. Then problem \eqref{Eq
sem. grp.} possesses the unique solution $u(t)$ such that
\begin{equation}\label{res 2}
u\in C([0,T]; D(A)) , \; u_t \in L^{\infty}(0,T; L^2(\Omega)),\ \
T>0.
\end{equation}
Moreover, if $u_0\in L^2(\Omega)$ and $f \in
L^1_{loc}(\R^+;L^2(\Omega)),$ then \eqref{Eq sem. grp.} possesses a
unique (mild) solution $u \in C([0,T]; L^2(\Omega)) $ given by
\begin{equation}\label{res 3}
u(t)=S(t)u_0 + \int_0^t S(t-s)f(s)\, ds
\end{equation}
where $S(t)$ is a semigroup of contractions generated by $A.$
\end{proposition}

\begin{corollary}\label{res. 2 teman} Under the hypothesys of Proposition \ref{res. 1 teman}, the solution
$u$ in (\ref{res 2}) satisfies
\begin{equation}\label{res 2.1}
u\in L^{\infty}\left( 0,T;H^1_0(\Omega)\cap H^2(\Omega)\right),
\end{equation}
\end{corollary}
For the proof, see \cite{temam2}.

Furthermore, one can get (see \cite{kato}, for instance) the
estimate for strong solution (\ref{res 2}):
\begin{equation}\label{estimativa temam/kato 1}
\|u_t \| (t) \leq \|Au_0\| + \|f\|(0) + \|f_t\|_{L^1_tL^2_{xy}},
\end{equation}
and
\begin{equation}\label{estimativa temam/kato 2}
\left\| A u \right\| (t) \leq \|  u_t \|(t)+ \|f\|(t).
\end{equation}
Since $D(A){ \hookrightarrow }  H^1_0(\Omega)\cap H^2(\Omega)$
compactly (see \cite{temam2} for instance), we have the estimate
\begin{eqnarray}\label{estimativa temam/kato 4}
\|u\|_{L^{\infty}\left(0,T; H^1_0\cap H^2 (\Omega)\right)}(t) \leq C
\big(\|u\|_{L^{\infty}_t L^2_{xy}} +  \|Au_0\| + \|f\|(0) +
\|f_t\|_{L^1_tL^2_{xy}}  +  \|f\| _{L^{\infty}_t L^2_{xy}} \big).
\end{eqnarray}
%=======================
where $C$ depends only on $\Omega $. Next, we define
$$
Y_T = \{f \in L^{1}\big( 0,T; L^2(\Omega) \big) \; \text{such that}
\; f_t \in L^{1} \big( 0,T; L^2(\Omega)   \big)     \}
$$
with the norm
$$
\|f\|_{Y_T}= \|f\|_{ L^{1}_t L^2_{xy} }  + \|f_t\|_{L^{1}_t L^2_{xy}
}.
$$
\begin{remark}\label{prop traÃ§o f} If  $f\in Y_T,$ then $ f \in C([0,T];L^2(\Omega))$
and the following inequality holds:
$$\sup_{t\in[0,T]}\|f\|_{L^2_{xy}} \leq C_T
\|f\|_{Y_T}
$$
where $C_T$ is proportional to $T$ and its positive powers
\cite{evans}.
\end{remark}
Consider
$$ X_T^0 =  \left\{ u \in L^{\infty}\big(0,T; H^1_0(\Omega)   \big)\, ; \,
\nabla u_y \in L^{\infty}\big(0,T; L^2(\Omega)\big)\, \text{ and
}\,u_{xx}\in L^2\big(0,T; L^2(\Omega)\big) \right\}$$ with the norm
$$
\|u\|_{X_T^0}= \|u\|_{ L^{\infty}_t H^1_{0xy} }  + \|\nabla u_y\|_{ L^{\infty}_t L^2_{xy} } + \|u_{xx}\|_{ L^2_t L^2_{xy} }
.
$$ and define the Banach space
\begin{eqnarray}\label{esp XT}
X_T=\{u \in X_T^0:\ u_t \in L^{\infty}\big(0,T; L^2(\Omega)\big)
\;\;\text{and}\;\; \nabla u_t \in L^{2}\big(0,T; L^2(\Omega)\big)
\} .
\end{eqnarray}
%-------------------------------
with the norm
\begin{eqnarray}\label{morma esp XT}
\| u \|_{X_T} = \|u\|_{X_T^0}
+ \| u_t \|_{L^{\infty}_T L^2_{xy}}
+\| \nabla u_t \|_{L^{2}_T L^2_{xy}} .
\end{eqnarray}
%-------------------------------
%Observe que  $L^{\infty}\big(0,T; H^1_0(\Omega) \cap H^2(\Omega)   \big)\hookrightarrow  X_T^0 ,$ com imersão dada por
%\begin{eqnarray}\label{morma esp XT}
%\| u \|_{X_T^0} \leq  \|u\|_{ L^{\infty}_t H^1_0\cap H^2{xy} }  .
%\end{eqnarray}
%-------------------------------

\begin{thm}\label{Teo solu local mZK} Let $\delta=1$ and $u_0 \in D(A)$. Then there exists $T>0$
such that IBVP (\ref{2.1})-(\ref{2.4}) possesses a unique
%(weak/mild)
solution in $X_T$.
\end{thm}
The proof of the Theorem consists in three lemmas below.
\begin{lemma}\label{lemma1} The function
%\begin{eqnarray*}
$Y_T
%&
\longrightarrow
%&
X_T;\ \
%\\
%-------------------------------
f
%&
\mapsto
%&
\int_0^t S(t-s)f(s) ds$
%\end{eqnarray*}
is well defined and
continuous.
\end{lemma}
%-------------------------------
%-------------------------------
%\noindent{{\sl Proof:}}
For the proof, note that this function maps $f$ to the solution of
homogeneous linear problem with zero initial datum. Estimates
(\ref{estimativa temam/kato 1}) and (\ref{estimativa temam/kato 4})
then give
\begin{equation}\label{estimativa ut 0}
\| u \|_{L^{\infty}_T H^1_0 \cap H^2_{xy} } + \| u_t
\|_{L^{\infty}_T L^2_{xy}} \leq C_1 \|f\|_{Y_T}.
\end{equation}
Hence,
\begin{equation}\label{estimativa ut 0.1}
\| u \|_{X_T^0} + \| u_t \|_{L^{\infty}_T L^2_{xy}} \leq C_2 \|f\|_{Y_T}.
\end{equation}
Thus, it rests to estimate the term $\|
\nabla u_t \|_{L^{2}_T L^2_{xy}}$ in (\ref{morma esp XT}).

Differentiate the equation in (\ref{Eq sem. grp.}) with respect to
$t,$ multiply it by $(1+x) u_{t}$ and integrate the outcome over
$\Omega.$ The result reads
\begin{equation}\label{estimativa ut}
\dfrac{d}{dt} \left( (1+x),u_t^2\right) (t) + \|\nabla u_t \|^2(t) +
2\| u_{xt} \|^2+\int_{-B}^B u_{xt}^2(0,y,t)\,dy = \|u_t\|^2(t) +
2\int_{\Omega}(1+x) f_t  u_{t} \, d\Omega.
\end{equation}
H\"older's inequality and (\ref{estimativa temam/kato 1}) imply
\begin{eqnarray}\label{estimativa ut 1.2}
\int_0^T  \|\nabla u_t \|^2(t)\, dt &\leq &  T\big(\|f\|(0) +
\|f_t\|_{L^1_T L^2_{xy}}  \big)^2 \nonumber \\
&+& 2(1+L)\big(\|f\|(0) + \|f_t\|_{L^1_T L^2_{xy}} \big)
\|f_t\|_{L^1_T L^2_{xy}} +\left( (1+x),u_t^2\right) (0) .
%==============================
\end{eqnarray}
%-------------------------------
%==============================
Using the equation from (\ref{Eq sem. grp.}) and taking in mind that
$u_0 \equiv 0 $, we get
\begin{eqnarray}\label{estimativa ut 1.3}
u_t (x,y,0) = f(x,y,0) - Au_0= f(x,y,0)
\end{eqnarray}
%=======================
Inserting (\ref{estimativa ut 1.3}) into (\ref{estimativa ut 1.2})
provides
\begin{eqnarray}\label{estimativa ut 1.7}
\|\nabla u_t \|^2_{L^2_TL^2_{xy}}  \leq  \Big(4T K_T^2 + 4K_T(1+L) + K_T^2(1+L)  \Big)\|f\|^2_{Y_T},
\end{eqnarray}
%==============================
where $K_T=\max \{1,C_T\}$. Therefore, estimates (\ref{estimativa ut
0}) and (\ref{estimativa ut 1.7}) read
\begin{eqnarray}\label{estimativa ut 1.9}
\|u\|_{X_T} \leq K \|f\|_{Y_T}.
\end{eqnarray}
%==============================
\begin{lemma} The function
$$
D(A)   \longrightarrow X_T ;\ u_0  \mapsto  S(t)u_0
$$
is well defined and continuous.
\end{lemma}
%-------------------------------
The proof follows the same steps as Lemma \ref{lemma1}, taking into
account that now $f\equiv 0$. The resulting estimate is
\begin{eqnarray}\label{estimativa linear XT}
 \|u\|_{X_T} & \leq &  M \|u_0 \|_{D(A)},
\end{eqnarray}
%-------------------------------
where $M$ is given by
\begin{eqnarray}\label{estimativa linear M}
M = 2C + 1 + \sqrt{1+L+T},
\end{eqnarray}
and $C$ (which depends only on $\Omega$) is defined by continuous
immersion $D(A) \hookrightarrow  H^1_0(\Omega)\cap H^2(\Omega).$
%-----------------------------------------------------------------

\begin{lemma}\label{contraÃ§Ã£o}
Given $R>0$, consider the closed ball $ B_R = \{u \in X_T  ;
\|u\|_{X_T} \leq R \}.$ If $R>0$ is sufficiently small, then the
operator
\begin{eqnarray}\label{contraÃ§ao 0}
\Phi : B_R  \longrightarrow  X_T ;\
%-------------------------------
v \mapsto S(t)u_0 -\int_0^t S(t-s)v^2 v_x (s) \,ds \nonumber
\end{eqnarray}
is the contraction.
\end{lemma}
%-------------------------------
Fix $R>0$ and $u,v \in B_R.$ We have
\begin{eqnarray*}
\Phi (v) - \Phi (u)= \int_0^t S(t-s)[u^2u_x - v^2 v_x] (s) \,ds
\nonumber
\end{eqnarray*}
%-------------------------------
so that (\ref{estimativa ut 1.9}) implies
\begin{eqnarray}\label{contraÃ§ao 0.1}
\| \Phi (u) - \Phi (v) \|_{X_T} \leq K \|u^2u_x - v^2v_x\|_{Y_T}.
\end{eqnarray}
%-------------------------------
We study the right-hand norm in detail:
\begin{eqnarray}\label{contaÃ§ao 1}
\|u^2u_x - v^2v_x\|_{Y_T} &=& \|u^2u_x - v^2v_x\|_{L^{1}_T L^2_{xy}}
+\left\|\big(u^2u_x\big)_t - \big(v^2v_x\big)_t\right\|_{L^{1}_Y L^2_{xy}} \nonumber \\
&=& I + J.
\end{eqnarray}
%=============================
First, we write
\begin{eqnarray}\label{contaÃ§ao 2}
I &=& \left\|(u^2 - v^2)u_x\right\|_{L^{1}_T L^2_{xy}} +\left\| v^2(u_x-v_x)\right\|_{L^{1}_T L^2_{xy}}   \nonumber \\
&=& I_1 + I_2 .
\end{eqnarray}
For the integral $I_1$ one has
\begin{eqnarray}\label{contaÃ§ao 2.1}
I_1 \leq  \int_0^T  \| u-v\|_{L^{6}(\Omega)} \|
u+v\|_{L^{6}(\Omega)} \|u_x \|_{L^{6}(\Omega)} dt.
\end{eqnarray}
%============================
Nirenberg's inequality gives
\begin{eqnarray}\label{contaÃ§ao 2.4}
I_1  &\leq & T^{\frac{2}{3}} C_{\Omega} \| \nabla (u  + v)\|^{\frac{2}{3}}_{L^{\infty}_T L^2_{xy}} \|u
+ v\|^{\frac{1}{3}}_{L^{\infty}_T L^2_{xy}} \| \nabla u_x\|^{\frac{2}{3}}_{L^{2}_T L^2_{xy}}
\|u_x\|^{\frac{1}{3}}_{L^{\infty}_T L^2_{xy}} \| \nabla (u - v)\|^{\frac{2}{3}}_{L^{\infty}_T L^2_{xy}}
\|u - v\|^{\frac{1}{3}}_{L^{\infty}_T L^2_{xy}} \nonumber \\
%=============================
&\leq & T^{\frac{2}{3}} C_{\Omega}  D^{\frac{2}{3}} \| u  + v\|_{X_T} \|u\|_{X_T} \|
u - v \|_{X_T} ,
\end{eqnarray}
%=============================
where $D$ is the Poincare's constant from $\|w\|\leq D \| \nabla
w\|.$ Since $u$ and $v$ lie in $B_R,$ we conclude
\begin{eqnarray}\label{contaÃ§ao 2.5}
I_1 \leq  T^{\frac{2}{3}}K_0R^2 \| u - v \|_{X_T}.
\end{eqnarray}
%=============================
The integral $I_2$ can be treated in the similar way as $I_1$. It
rests to estimate the integral $J$.
\begin{eqnarray}\label{contaÃ§ao 3}
J \leq  \| 2u u_t(u_x - v_x)  \|_{L^1_T L^2_{xy}} + \| u^2 (u_{xt} - v_{xt}) \|_{L^1_T L^2_{xy}} + \|2v_x u(u_t - v_t)  \|_{L^1_T L^2_{xy}} \nonumber \\
%=============================
+ \|2v_x v_t (u-v)  \|_{L^1_T L^2_{xy}} +\| v_{xt}(u-v)(u+v) \|_{L^1_T L^2_{xy}} \nonumber \\
%=============================
= J_1 + J_2 + J_3 +J_4 + J_5.
\end{eqnarray}
%=============================
For $J_1$ we have
\begin{eqnarray}\label{contaÃ§ao 3.1}
J_1 & \leq \int_0^T  \| u\|_{L^{6}(\Omega)} \| u_t\|_{L^{6}(\Omega)} \|u_x - v_x \|_{L^{6}(\Omega)} \, dt .
\end{eqnarray}
%=============================
Niremberg's inequality implies
\begin{eqnarray}\label{contaÃ§ao 3.2}
J_1  \leq T^{\frac{1}{3}}   C_{\Omega} \| \nabla u\|^{\frac{2}{3}}_{L^{\infty}_T L^2_{xy}} \|u\|^{\frac{1}{3}}_{L^{\infty}_T L^2_{xy}} \| \nabla (u_x - v_x)\|^{\frac{2}{3}}_{L^{2}_T L^2_{xy}} \|u_x-v_x\|^{\frac{1}{3}}_{L^{\infty}_T L^2_{xy}}\|u_t\|^{\frac{1}{3}}_{L^{\infty}_T L^2_{xy}}\|\nabla u_t\|^{\frac{2}{3}}_{L^{2}_T L^2_{xy}} \nonumber \\
\leq T^{\frac{1}{3}} K_2 R^2  \| u - v \|_{X_T}.
\end{eqnarray}
%=============================
The integrals $J_3$ and $J_4$ are analogous to $J_1$. To get bound
for $J_5$ we observe that
\begin{eqnarray}\label{contaÃ§ao 3.6}
J_5  &=& \int_0^T \left(\int_{\Omega} v_{xt}^2 (u -v)^2 (u+v)^2 \, d\Omega\right)^{\frac{1}{2}} dt \nonumber \\
%==============================
& \leq & \int_0^T \left( \sup_{(x,y) \in \Omega} (u-v)^2(x,y,t) \right)^{\frac{1}{2}} \left( \sup_{(x,y) \in \Omega} (u+ v)^2(x,y,t) \right)^{\frac{1}{2}} \|v_{xt}\|(t) \, dt \nonumber \\
%==============================
&\leq & \int_0^T \big( \|u-v\|^2_{H^1_{xy}}(t) + \|u_{xy} - v_{xy}\|^2 (t) \big)^{\frac{1}{2}} \big( \|u+v\|^2_{H^1_{xy}}(t) + \|u_{xy} + v_{xy}\|^2 (t) \big)^{\frac{1}{2}}   \| v_{xt}\|(t) \, dt \nonumber \\
%=============================
&\leq & \big( \|u-v\|_{L^{\infty}_T H^1_{xy}} + \|u_{xy}- v_{xy}\|_{L^{\infty}_T L^2_{xy}} \big)\big( \|u+v\|_{L^{\infty}_T H^1_{xy}} + \|u_{xy}+ v_{xy}\|_{L^{\infty}_T L^2_{xy}} \big)\| v_{xt}\|_{L^{1}_T L^2_{xy}} \nonumber \\
%=============================
& \leq &  T^{\frac{1}{2}}C_{\Omega} \| v\|_{X_T}\|u+v\|_{X_T} \|u-v\|_{X_T} \nonumber \\
%=============================
& \leq & 2 T^{\frac{1}{2}} C_{\Omega}  R^2 \|u-v\|_{X_T}.
\end{eqnarray}
%=============================
The integral $J_2$ follows like $J_5$. Thus,
\begin{eqnarray}\label{contaÃ§ao 3.6}
\|u^2u_x - v^2v_x\|_{Y_T} \leq K K^* T^{\frac{1}{3}}R^2 \|u-v\|_{X_T}.
\end{eqnarray}
%=============================
Finally, choosing $T>0$ such that $K K^* T^{\frac{1}{3}}R^2<1 ,$ we
conclude that $\Phi $ is a contraction map.

Lemma \ref{contraÃ§Ã£o} is proved.

Let $u \in B_R$. If $R= 2 M\|u_0\|_{D(A)},$ then estimates
(\ref{estimativa linear XT}) and (\ref{contaÃ§ao 3.6}) with $v
\equiv 0$ assure
\begin{eqnarray}\label{contaÃ§ao 3.7}
\|u\|_{X_T} & \leq & \|S(t)u_0\|_{X_T} +\| \int_0^t S(t-s)u^2u_x \, ds \|_{X_T} \nonumber \\
%=============================
& \leq & M\|u_0\|_{D(A)} + K K^* T^{\frac{1}{3}}R^2 \|u\|_{X_T} \nonumber \\
%=============================
&\leq & \frac{R}{2} + K K^* T^{\frac{1}{3}}R^3.
\end{eqnarray}
%=============================
Setting $T>0$ such that $K K^* T^{\frac{1}{3}}R^3 < \frac{R}{2}, $
one get
\begin{eqnarray}\label{contaÃ§ao 3.8}
\|u\|_{X_T} \leq R.
\end{eqnarray}
Choose $T>0$ such that $K K^* T^{\frac{1}{3}}R^2<1 $ and $ K K^*
T^{\frac{1}{2}}R^3 < \frac{R}{2} .$ Then $\Phi$ is the contraction
from the ball $B_R$ into itself. Therefore, the Banach fixed point
theorem assures the existence of a unique element $u\in B_R$ such
that $\Phi (u) = u.$

This completes the proof of Theorem \ref{Teo solu local mZK}.

\section{Global estimates and decay}\label{global}
%-------------------------------

\begin{thm}\label{teorema principal}
Let $B,L>0$ and $u_0(x,y)$ be such that
$$\frac{2\pi^2}{L^2}-1 >0, \ \ A^2:=\frac{\pi^2}2
\left[\frac{3}{L^2}+\frac{1}{4B^2}\right] - 1 >0 \ \ \text{ and
 }\ \ \|u_0\|^2 < \dfrac{A^2}{2\pi^2\left(
\frac{1}{L^2}+\frac{1}{4B^2}\right)}.$$

Suppose $u_0\in D(A)$ satisfies
$$
I_0^2=\|u_{0x} + \Delta u_{0x} + u_0^2u_{0x}\|^2 < \infty,
$$
and
\begin{eqnarray}\label{cond sobre u0}
 \left[ \frac{2(1+L)^2}{1-2\|u_0\|^2}\|u_0\|^2 \big( I_0^2 + \|u_0\|^2  \big)\right]
 \left[4^2 + \frac{6^3(4!)^2(1+L)^8}{(1-2\|u_0\|^2)^2}\big( I_0^2 + \|u_0\|^2  \big)^2 \right]
 < \frac{2\pi^2}{L^2}-1.
\end{eqnarray}
Then for all $T>0$ there exists a unique solution $u\in X_T$ to
problem \eqref{2.1}-\eqref{2.4}; more precisely,
%\begin{eqnarray*}
$$
u\in L^{\infty}\left(0,T;H^1_0(\Omega)\right),  \ \nabla u_y \in L^{\infty}\left(0,T;L^2\Omega)\right),
\ u_{xx} \in L^2\left(0,T;L^2(\Omega)\right),
%\\
$$$$u_t \in L^{\infty}\left(0,T;L^2\Omega)\right), \ \nabla u_t \in
L^2\left(0,T;L^2(\Omega)\right). $$
%\end{eqnarray*}
Moreover, there exist constants $C>0$ and $\gamma >0$ such that
\begin{equation}\label{decaimento 1}
%-------------------------------
\|u\|^2_{H^1(\Omega)}(t) + \|\nabla u_y\|^2(t)+ \| u_t\|^2(t) \leq C
e^{-\gamma t},\ \ \forall t\ge 0
%-------------------------------
\end{equation}
and, in addition,
$$
u_x(0,y,t),\ u_{xy}(0,y,t),\ u_{xx}(L,y,t) \in
L^{\infty}\left(0,T;L^2(-B,B)\right), $$$$
 u_{xx}(0,y,t)\in L^2\left(0,T;L^2(-B,B)\right).
$$
\end{thm}

%-------------------------------
%-------------------------------

Let $u\in X_{T_0}$ be a local solution given by Theorem (\ref{Teo
solu local mZK}). We are going to obtain a priori estimates
independent of $T_{0}$ in order to extend the solution to all $T>0.$

\subsection{Estimate I}
We start the proof of \eqref{decaimento 1}, multiplying \eqref{2.1}
by $u$ and integrating over $Q_t,$ which easily gives

\begin{equation}\label{decaimento 1.3}
\| u\| ^2 (t) \leq \|u_0\|^2.
\end{equation}

\subsection{Estimate II}
Multiplying \eqref{2.1} by $(1+x)u$ and integrating over $\Omega,$
we have

\begin{eqnarray}\label{decaimento 1.6}
\dfrac{d}{dt} \left( 1+x, u^2 \right)(t) + \int_{-B}^B u^2_x (0,y,t)
\, dy
+  \|\nabla u \|^2(t) + 2\|u_x\|^2(t) - \|u\|^2(t) \nonumber \\
%-------------------------------
=   -2\int_{\Omega}(1+ x) u (u^2u_x) \, d\Omega %\nonumber \\
%-------------------------------
=  \frac{1}{2}\int_{\Omega}u^4 \, d\Omega.%\nonumber \\
\end{eqnarray}
%-------------------------------
%-------------------------------
For the integral $I_1 = \frac{1}{2}\int_{\Omega}u^4 =
\frac{1}{2}\|u\|^4_{L^4(\Omega)}(t),$
%\begin{eqnarray}\label{decaimento 1.7}
%\end{eqnarray}
%-------------------------------
%-------------------------------
Nirenberg's inequality yields

\begin{eqnarray}\label{decaimento 1.8}
I_1 &\leq & \frac{1}{2} \big(  2^{\frac{1}{2}} \|\nabla u \|^{\frac{1}{2}}(t)\| u \|^{\frac{1}{2}} (t) \big)^4 \nonumber \\
%-------------------------------
&=& 2 \|\nabla u \|^{2}(t)\| u \|^{2}(t)%\nonumber \\
%-------------------------------
\leq  2 \|\nabla u \|^{2}(t)\| u_0\|^{2}(t).
\end{eqnarray}

Take
$$ I_2 = 3\|u_x\|^2(t)+ \|u_y\|^2(t).
$$
For all $\varepsilon>0$ we have
$$ I_2 = (3-\varepsilon)\|u_x\|^2(t)+ (1-\varepsilon)\|u_y\|^2(t) + \varepsilon\big( \|u_x\|^2(t)+ \|u_y\|^2(t) \big).
$$
%-------------------------------
%-------------------------------
Lemma \ref{prop 1 decaimento} jointly with (\ref{decaimento 1.6})
and (\ref{decaimento 1.8}) provides

\begin{eqnarray}\label{decaimento 1.9}
\dfrac{d}{dt} \left( 1+x, u^2 \right)(t)
+ \left[ \pi^2\left( \frac{3}{L^2}+\frac{1}{4B^2} \right) -1 - \varepsilon \pi^2\left( \frac{1}{L^2}+\frac{1}{4B^2} \right)  \right]\|u\|^2(t)  \nonumber \\
%-------------------------------
+ \left( \varepsilon - 2\|u_0\|^2 \right)\|\nabla u\|^2(t) \leq 0.
\end{eqnarray}
%-------------------------------
%-------------------------------
Define
$$
2 A^2:=\pi^2 \left[\frac{3}{L^2}+\frac{1}{4B^2}\right] - 1 >0,
%$$
\ \text{ and take }\
%$$
\varepsilon= \dfrac{A^2}{\pi^2\left(
\frac{1}{L^2}+\frac{1}{4B^2}\right)}.
$$
The result for \eqref{decaimento 1.9} reads
\begin{eqnarray}\label{decaimento 1.10}
\dfrac{d}{dt} \left( 1+x, u^2 \right)(t) + A^2\|u\|^2(t)+ \left(
\varepsilon - 2\|u_0\|^2 \right)\|\nabla u\|^2(t) \leq 0.
\end{eqnarray}
%-------------------------------
%-------------------------------
If $0\leq \varepsilon - 2\|u_0\|^2,$ then
\begin{eqnarray}\label{decaimento 1.11}
\dfrac{d}{dt} \left( 1+x, u^2 \right)(t) + \frac{A^2}{1+L}\left(
1+x, u^2 \right)(t) \leq 0,
\end{eqnarray}
%-------------------------------
%-------------------------------
and consequently
\begin{equation}\label{decaimento 1.12}
%-------------------------------
\|u\|^2(t) \leq \left( 1+x, u^2 \right) (t) \leq e^{-\gamma_0
t}\left( 1+x, u^2_0 \right)\text{  with }\gamma_0 =\frac{A^2}{1+L}.
%-------------------------------
\end{equation}

\subsection{Estimate III}
Write (\ref{decaimento 1.6}) as
\begin{eqnarray}\label{decaimento 4}
\|\nabla u \|^2(t) + 2\|u_x\|^2(t)+\int_{-B}^B u^2_x (0,y,t) \, dy  = -2
\left( (1+x)u,u_t + u_x + u^2u_x\right) (t) \nonumber \\
%=================
=  -2\left( (1+x)u,u_t + u_x \right) (t) +\frac{1}{2}\|u\|^4_{L^4_{xy}}(t)\nonumber \\
%=================
\leq 2(1+L)\|u\|(t)\big( \|u_t\|(t) + \|u_x\|(t)  \big) + 2\|u\|^2(t)\|\nabla u\|^2(t).
\end{eqnarray}
%========================
Then
\begin{eqnarray}\label{decaimento 4.1}
(1-2\|u_0\|^2)\|\nabla u \|^2(t) + \|u_x\|^2(t)+\int_{-B}^B u^2_x
(0,y,t) \, dy \leq 2(1+L)^2\big( \|u_t\|^2(t) + \|u\|^2(t)  \big).
\end{eqnarray}

Note for posterior use that $\nabla u$ is estimated by $u$ e $u_t$
provided $u_0$ be sufficiently small in $L^2(\Omega)$:

\begin{eqnarray}\label{decaimento 4.2}
\|\nabla u \|^2(t)  \leq C_{\|u_0\|} \big( \|u_t\|^2(t) + \|u\|^2(t)
\big)
\end{eqnarray}
%========================
where
\begin{eqnarray}\label{decaimento 4.3}
 C_{\|u_0\|}=\frac{2(1+L)^2}{(1-2\|u_0\|^2) }.
\end{eqnarray}

\subsection{Estimate IV} Differentiate the equation with respect to $t$ and multiply the result by $(1+x) u_{t}$.
Integrating over $\Omega$ then gives

\begin{equation}\label{decaimento ut}
\dfrac{d}{dt} \left( (1+x),u_t^2\right) (t) + \|\nabla u_t \|^2(t) +
2\| u_{xt} \|^2+\int_{-B}^B u_{xt}^2(0,y,t)\,dy =
\|u_t\|^2(t)-\frac{2}{3}\int_{\Omega}(1+x) (u^3)_{xt}  u_{t} \,
d\Omega.
\end{equation}
We have
\begin{equation}\label{decaimento ut 1}
-\frac{2}{3}\int_{\Omega}(1+x) (u^3)_{xt}  u_{t} \, d\Omega = \left(
u^2,u_t^2 \right)(t) - 2\left( (1+x)uu_x, u_t^2 \right) (t) = I_1 +
I_2.
\end{equation}
H\"{o}lder and  Nirenberg's inequalities provide
\begin{eqnarray}\label{decaimento ut 2}
 I_1 \leq \|u\|^2_{L^4_{xy}} \|u_t\|^2_{L^4_{xy}}
 &\leq & 4 \|u\|(t)\|\nabla u\|(t)\|u_t\|(t)\|\nabla u_t\|(t)\nonumber \\
& \leq & \frac{1}{4}\|\nabla u_t\|^2(t) + 4^2\|u\|^2(t)\|\nabla
u\|^2(t)\|u_t\|^2(t).
\end{eqnarray}
Nirenberg's inequality for $p=8$ then implies
\begin{eqnarray}\label{decaimento ut 3}
 I_2 &\leq & 2(1+L)\|u_x\|(t)\|u\|_{L^4_{xy}}(t)\|u_t\|_{L^8_{xy}}^2(t)\nonumber \\
%==================
 &\leq & 2^{\frac{3}{2}}C_{N8}^2(1+L)\|u_x\|(t)\|u\|^{\frac{1}{2}}(t)\|
 \nabla u\|^{\frac{1}{2}}(t)\|u_t\|^{\frac{1}{2}}(t)\|\nabla u_t\|^{\frac{3}{2}}(t)\nonumber \\
%====================
&\leq & 2^{\frac{3}{2}}C_{N8}^2(1+L)\|u\|^{\frac{1}{2}}(t)\|\nabla u\|^{\frac{3}{2}}(t)
\|u_t\|^{\frac{1}{2}}(t)\|\nabla u_t\|^{\frac{3}{2}}(t).\nonumber \\
\end{eqnarray}
Taking $k=\frac{4}{3}$ requires $l=4$ and by Young's inequality this
reads
\begin{eqnarray}\label{decaimento ut 4}
 I_2 \leq  \frac{1}{4}\|\nabla u_t\|^2(t) + C_1\|u\|^2(t)\|\nabla u\|^6(t)\|u_t\|^2(t)
\end{eqnarray}
where
$$
C_1= 2^4 3^3 C_{N8}^8(1+L)^4,
$$
and by Corollary \ref{corollary Nirenberg} we obtain that
\begin{eqnarray}\label{decaimento ut 5}
C_1= 2\cdot 3^3 (4!)^2(1+L)^4.
\end{eqnarray}
Thus,
\begin{eqnarray}\label{decaimento ut 6}
I_1 + I_2 \leq  \frac{1}{2}\|\nabla u_t\|^2(t) +\Big[
\|u\|^2(t)\|\nabla u\|^2(t)\Big]
\Big[4^2 + C_1 \|\nabla u\|^4(t)
\Big]\|u_t\|^2(t).
\end{eqnarray}
Using (\ref{decaimento 4.2}) we have
\begin{eqnarray}\label{decaimento ut 7}
I_1 + I_2 \leq  \frac{1}{2}\|\nabla u_t\|^2(t) +\Big[ C_{\|u_0\|}\|u\|^2(t) \big( \|u_t\|^2(t) + \|u\|^2(t)
 \big)\Big]\nonumber \\
%=================
\Big[4^2 + C_1 C_{\|u_0\|}^2\big( \|u_t\|^2(t) + \|u\|^2(t)  \big)^2
\Big]\|u_t\|^2(t).
\end{eqnarray}
Backing to (\ref{decaimento ut}), we get
\begin{eqnarray}\label{decaimento ut 8}
\dfrac{d}{dt} \left( (1+x),u_t^2\right) (t) +\frac{1}{2} \|\nabla
u_t \|^2(t) + 2\| u_{xt} \|^2+\int_{-B}^B u_{xt}^2(0,y,t)\,dy -
\big[1 + \omega(t)\big]\|u_t\|^2(t) \leq 0
\end{eqnarray}
with
\begin{eqnarray}\label{decaimento ut 9}
\omega(t)= \Big[ C_{\|u_0\|}\|u\|^2(t) \big( \|u_t\|^2(t) +
\|u\|^2(t)  \big)\Big]\Big[4^2 + C_1 C_{\|u_0\|}^2\big( \|u_t\|^2(t)
+ \|u\|^2(t)  \big)^2 \Big].
\end{eqnarray}
The use of Steklov's inequality gives
\begin{eqnarray}\label{decaimento ut 10}
\dfrac{d}{dt} \left( (1+x),u_t^2\right) (t) +\frac{\pi^2}{2(1+L)}\left(\frac{1}{L^2} +\frac{1}{4B^2} \right)
\left( 1+x,  u_t^2\right) (t)+\int_{-B}^B u_{xt}^2(0,y,t)\,dy \nonumber \\
%==================
+ \frac{1}{(1+L)}\left[\frac{2\pi}{L^2} -1 - \omega(t)\right]\left(
(1+x),u_t^2\right) (t) \leq 0 .
\end{eqnarray}

Setting $z(t)= \left( (1+x),u_t^2\right) (t),$ \eqref{decaimento ut
10} reads

\begin{eqnarray}\label{estb. Lyap. z(t)}
\dfrac{d}{dt} z(t) & \leq & \frac{1}{1+L}\left[1- \frac{\pi^2}{2}\left(\frac{5}{L^2} +\frac{1}{4B^2} \right)  \right] z(t)+  \omega(t)z(t) \nonumber \\
%======================
&=& p_1 z(t) + p_2 z^2(t) + p_3 z^3(t) + p_4 z^4(t)
\end{eqnarray}
where
 \begin{eqnarray}\label{estb. Lyap. z(t) 1}
p_1 &=&  \frac{1}{1+L}\left[1- \frac{\pi^2}{2}\left(\frac{5}{L^2} +\frac{1}{4B^2} \right)  \right] + C_{\|u_0\|}\|u_0\|^4 \left(4^2+C_1 C_{\|u_0\|}^2\|u_0\|^4 \right) \\
%=====================
p_2 &=&  C_{\|u_0\|}\|u_0\|^4 \left( 4^2+3 C_1 C_{\|u_0\|}^2\|u_0\|^2\right)       \\
%=====================
p_3 &=&  C_1C_{\|u_0\|}^3\|u_0\|^4   \left( 1+2\|u_0\|^2\right)      \\
%=====================
p_4 &=&C_1C_{\|u_0\|}^3\|u_0\|^2.
\end{eqnarray}

Next we compute $\frac{d}{dt}z(0)$ to show that $z(t)$ reaches a
local (lateral) maximum at $t=0.$ In order $\frac{d}{dt}z(0)$ to be
negative, it should be
%$p_1$ given by (\ref{estb. Lyap. z(t) 1}) be negative
\begin{equation}\label{estb. Lyap. z(t) 4}
 C_{\|u_0\|}\|u_0\|^4 \left(4^2+C_1 C_{\|u_0\|}^2\|u_0\|^4 \right)
 <  \frac{1}{1+L}\left[ \frac{\pi^2}{2}\left(\frac{5}{L^2} +\frac{1}{4B^2} \right) -1
 \right].
\end{equation}

%Since $u_0$ is arbitrary,
Without loss of generality
one can assume
%\begin{eqnarray}\label{decaimento ut 11}
$\left((1+x),u_t^2\right) (0)  > 0.$
%\end{eqnarray}
Choose $u_0 \in D(A)$ such that \eqref{cond sobre u0} holds, i.e., $
\omega(0) <\frac{2\pi^2}{L^2}-1.$ Then
$$p_1 (1+L)I_0^2 + p_2 [(1+L)I_0^2]^2 + p_3 [(1+L)I_0^2]^3 + p_4 [(1+L)I_0^2]^4 <0,$$
which assures
\begin{eqnarray}\label{decaimento ut 16}
\frac{d}{dt}\left( (1+x),u_t^2\right) (0)  < 0.
\end{eqnarray}
This means that $\left( (1+x),u_t^2\right) (0)$ is a local
(left-hand) straight maximum.

Observe that
\begin{eqnarray}\label{decaimento ut 17}
   \|u_t\|^2(t) \leq (1+L)\|u_t\|^2(0) , \,\,\,\,\,\,\,\,\,\forall t \in [0,T_0) .
\end{eqnarray}
Therefore,
\begin{eqnarray}\label{decaimento ut 18}
\omega(t) &= & \Big[ C_{\|u_0\|}\|u\|^2(t) \big( \|u_t\|^2(t) + \|u\|^2(t)  \big)\Big]
 \Big[4^2 + C_1 C_{\|u_0\|}^2\big( \|u_t\|^2(t) + \|u\|^2(t)  \big)^2 \Big]\nonumber \\
%===================
&\leq & \Big[ C_{\|u_0\|}\|u\|^2(0) \big( (1+L)\|u_t\|^2(0) + \|u\|^2(0)  \big)\Big]
 \Big[4^2 + C_1 C_{\|u_0\|}^2\big( (1+L)\|u_t\|^2(0) + \|u\|^2(0)  \big)^2 \Big]\nonumber \\
%====================
&\leq & (1+L)\Big[ C_{\|u_0\|}\|u\|^2(0) \big( \|u_t\|^2(0) + \|u\|^2(0)  \big)\Big]
 \Big[4^2 + C_1 C_{\|u_0\|}^2\big( \|u_t\|^2(0) + \|u\|^2(0)  \big)^2 \Big]\nonumber \\
%====================
&=& (1+L)\omega(0).
\end{eqnarray}
Consequently,

\begin{eqnarray}\label{decaimento ut 19}
 \frac{1}{(1+L)}\left[\frac{2\pi^2 }{L^2}-1  - \omega(t)\right]& = & \left[\frac{2\pi^2-L^2}{L^2(1+L)}
  -\frac{1}{(1+L)} \omega(t)\right] \nonumber \\
%===============
&\geq&   \left[\frac{2\pi^2-L^2}{L^2(1+L)}  -\omega(0)\right]\nonumber \\
&>&  0.
\,\,\,\,\,\,\,\,\,\,\,\,\,\,\,\,\,\,\,\,\,\,\,\,\,\,\,\,\,\,\,\,\,\,\,\,\,\,\,\,
\forall t \in [0,T_0).
\end{eqnarray}
Integrating the inequality (\ref{decaimento ut 10}) gives
\begin{eqnarray}\label{decaimento ut 19}
\left( (1+x),u_t^2\right) (T_0) < \left( (1+x),u_t^2\right) (0) ,
\end{eqnarray}
Then $\frac{d}{dt} \left( (1+x),u_t^2\right) (T_0)  < 0.$

Using (\ref{decaimento ut 18}) in (\ref{decaimento ut 10}), we have
\begin{eqnarray}\label{decaimento ut 20}
\dfrac{d}{dt} \left( (1+x),u_t^2\right) (t) +\frac{\pi^2}{2(1+L)}\left(\frac{1}{L^2} +\frac{1}{4B^2} \right)\left( 1+x,  u_t^2\right) (t)+\int_{-B}^B u_{xt}^2(0,y,t)\,dy \nonumber \\
%==================
+ \frac{1}{(1+L)}\left[\frac{2\pi^2}{L^2} -1 -
\omega(0)\right]\left( (1+x),u_t^2\right) (t) \leq 0 , \,\,\,\,\,\,
\forall t \in [0,T].
\end{eqnarray}

Finally,
\begin{eqnarray}\label{decaimento ut 21}
 \|u_t\|^2 (t)  \leq (1+L)\|u_t\|^2(0)e^{-\gamma_1 t}, \,\,\,\,\,\, \forall t \in [0,T]
\end{eqnarray}
where
\begin{eqnarray}\label{decaimento ut 22}
\gamma_1= \frac{\pi^2}{2(1+L)}\left(\frac{1}{L^2} +\frac{1}{4B^2} \right).
\end{eqnarray}
Integrating (\ref{decaimento ut 20}) over $[0,T],$ we obtain that
\begin{eqnarray}\label{decaimento ut 22.1}
\|\nabla u_t\|^2_{L^2_T L^2_{xy}} \leq \frac{2(1+L)^2}{\pi^2}\left(\frac{1}{L^2} +\frac{1}{4B^2} \right)^{-1}\|u_t\|^2(0)
\end{eqnarray}

Backing to(\ref{decaimento 4.2}) gives

\begin{eqnarray}\label{decaimento ut 23}
\|\nabla u \|^2(t)  &\leq & C_{\|u_0\|}  (1+L)\left(\|u_t\|^2(0)e^{-\gamma_1 t}
+ \| u_0 \|^2 e^{-\left(\frac{A^2}{(1+L)}\right)t} \right)\nonumber \\
%================
&\leq & C_{\|u_0\|}  (1+L)K_1 e^{-\gamma_2 t}, \,\,\,\,\,\, \forall t \in [0,T],
\end{eqnarray}
%========================
where

\begin{eqnarray}\label{decaimento ut 24}
\gamma_2 = \min \left\{\gamma_0,\gamma_1 \right\}
\end{eqnarray}
with $\gamma_0$ defined in \eqref{decaimento 1.12} and
\begin{eqnarray}\label{decaimento ut 25}
K_1 =\max \left\{\|u_t\|^2(0),\| u_0 \|^2 \right\}.
\end{eqnarray}

\subsection{Estimate V}
Multiply the equation by $(1+x) u_{yy}$ and integrate over $\Omega$.
The result reads

\begin{equation}\label{decaimento uyy}
 \|\nabla u_y \|^2(t) + 2\| u_{xy} \|^2+\int_{-B}^B u_{xy}^2(0,y,t)\,dy
 =2 \left( (1+x)u_{yy},u_t\right) (t) + \|u_y\|^2(t) + 2 \left( (1+x)u_{yy},u^2u_x\right) (t) .
\end{equation}
%============
We have
\begin{equation}
2 \left((1+x)u_{yy},u^2u_x\right) (t)  = \left( u^2 , u_y^2 \right)
(t) - 2 \left( (1+x)uu_x, u_y^2 \right) (t) = I_1 + I_2 .
\end{equation}
%==============

H\"{o}lder and Nirenberg's inequalities imply
\begin{eqnarray}\label{decaimento uyy 1}
 I_1 \leq \|u\|^2_{L^4_{xy}} \|u_y\|^2_{L^4_{xy}} &\leq & 2C_{\Omega}^2 \|u\|(t)\|\nabla u\|(t)\|u_y\|(t)\|\nabla u_y\|(t)\nonumber \\
%====================
& \leq & \frac{1}{6}\|\nabla u_y\|^2(t) + 6C_{\Omega}^4\|u\|^2(t)\|\nabla u\|^4(t)
\end{eqnarray}
In the same manner, the H\"{o}lder and  Nirenberg inequalities with
$p=8$ provide
\begin{eqnarray}\label{decaimento uyy 2}
 I_2 &\leq & 2(1+L)\|u_x\|(t)\|u\|_{L^4_{xy}}(t)\|u_y\|_{L^8_{xy}}^2(t)\nonumber \\
%==================
 &\leq & 2^{\frac{3}{2}}C_{\Omega_8}^2(1+L)\|u_x\|(t)\|u\|^{\frac{1}{2}}(t)\|\nabla u\|^{\frac{1}{2}}(t)\|u_y\|^{\frac{1}{2}}(t)\|\nabla u_y\|^{\frac{3}{2}}(t)\nonumber \\
%====================
&\leq & 2^{\frac{3}{2}}C_{\Omega_8}^2(1+L)\|u\|^{\frac{1}{2}}(t)\|\nabla u\|^2(t)\|\nabla u_y\|^{\frac{3}{2}}(t)\nonumber \\
\end{eqnarray}
Setting $k=\frac{4}{3}$ and $l=4$ and applying generalized Young's
inequality we come to
\begin{eqnarray}\label{decaimento uyy 3}
 I_2 \leq  \frac{1}{6}\|\nabla u_y\|^2(t) + C_2\|u\|^2(t)\|\nabla u\|^8(t)
\end{eqnarray}
where
\begin{eqnarray}\label{decaimento uyy 4}
C_2= 2\cdot 3^6 C_{\Omega_8}^8(1+L)^4.
\end{eqnarray}
Hence,
\begin{eqnarray}\label{decaimento uyy 5}
I_1 + I_2 \leq  \frac{1}{3}\|\nabla u_y\|^2(t) +\Big[ \|u\|^2(t)\|\nabla u\|^4(t)\Big]\cdot \Big[4C_{\Omega}^4 + C_2 \|\nabla u\|^4(t) \Big].
\end{eqnarray}

In turn, (\ref{decaimento uyy}) becomes
\begin{eqnarray}\label{decaimento uyy 6}
 \frac{1}{2}\|\nabla u_y\|^2(t) + 2\| u_{xy} \|^2+\int_{-B}^B u_{xy}^2(0,y,t)\,dy
 & \leq & 3(1+L)^2 \|u_t\|^2(t) + \|u_y\|^2(t)  \nonumber \\
 %======================
 & & + 4C_{\Omega}^4  \|u\|^2(t)\|\nabla u\|^4(t)+  C_2\|u\|^2(t) \|\nabla u\|^8(t)  \nonumber \\
 %======================
& \leq & (3(1+L)^2 + 1  )(1+L)K_1e^{-\gamma_2 t} \nonumber \\
 %======================
 & &+ 4C_{\Omega}^4  \|u_0\|^2   C_{\|u_0\|}^2  (1+L)^2K_1^2 e^{-2\gamma_2 t} \nonumber \\
 %======================
 & &+ C_2 \|u_0\|^2 C_{\|u_0\|}^4 (1+L)^4K_1^4 e^{-4\gamma_2 t}
  \nonumber \\
 %======================
 &\leq & K_2 e^{-\gamma_2 t}, \,\,\,\,\,\,\,\,\,\,\,\,\,\,\,\,\,\,\,\,\,\,\,\,\, \forall t \in [0,T]
\end{eqnarray}
%============
where
\begin{eqnarray}\label{decaimento uyy 7}
K_2 = \Big[ (3(1+L)^2 + 1  )(1+L)K_1 +  C_{\|u_0\|}^2  (1+L)^2K_1^2\|u_0\|^2
\left( 4C_{\Omega}^4   +C_2 C_{\|u_0\|}^2  (1+L)^2K_1^2 \right) \Big].
\end{eqnarray}

\subsection{Estimate VI} Now we have to estimate traces $u_x(0,y,t), u_{xy}(0,y,t)$ e $u_{xx}(L,y,t)$
in order to obtain the estimate for $\nabla u_x \in
L^{\infty}(0,T;L^2(\Omega)).$ From (\ref{decaimento 4.1}) we deduce
that
\begin{eqnarray}\label{decaimento 6}
\int_{-B}^B u^2_x (0,y,t) \, dy  &\leq &2(1+L)^2\big( \|u_t\|^2(t) + \|u\|^2(t)  \big) \nonumber \\
%==========================
&\leq &2(1+L)^3K_1e^{-\gamma_2 t}.
\end{eqnarray}
%========================
Now (\ref{decaimento uyy 6}) becomes
\begin{eqnarray}\label{decaimento 6.1}
\int_{-B}^B u_{xy}^2(0,y,t)\,dy  \leq K_2 e^{-\gamma_2 t}.
\end{eqnarray}
%============
Multiply it by $x$ and integrate over $(0,L).$ The result reads
\begin{eqnarray}\label{estimativa traço uxx L}
  Lu_{xx}(L,y,t) =  \|  u   + u_{yy} -x u_t -xu^2u_x\|_{L^1(0,L)} - u_x(0,y,t).
\end{eqnarray}
%-------------------------------
Therefore
\begin{eqnarray}\label{estimativa traço uxx L 1}
  L^2 u^2_{xx}(L,y,t)& \leq & 2\|  u   + u_{yy} -x u_t-xu^2u_x \|_{L^1(0,L)}^2 +2 u_x^2(0,y,t) \nonumber \\
  %================
 & \leq & 2L\|  u   + u_{yy} -x u_t -xu^2u_x\|_{L^2(0,L)}^2 +2 u_x^2(0,y,t) \nonumber \\
%==============================
 & \leq & 2L\Big( \|  u \|_{L^2(0,L)}  + \|u_{yy} \|_{L^2(0,L)}+\| x u_t \|_{L^2(0,L)} +\|xu^2u_x\|_{L^2(0,L)}\Big)^2 +2 u_x^2(0,y,t) \nonumber \\
%==============================
&= & 8 L \Big( \|  u \|_{L^2(0,L)}^2  + \|u_{yy} \|_{L^2(0,L)}^2+L^2\|  u_t \|_{L^2(0,L)}^2 +L^2\|u^2u_x\|_{L^2(0,L)}^2\Big)+2 u_x^2(0,y,t) \nonumber \\
%================================
\end{eqnarray}
%-------------------------------
Then

\begin{eqnarray}\label{estimativa traço uxx L 2}
  \int_{-B}^B  u^2_{xx}(L,y,t)\, dy  \leq  \frac{8}{L}  \Big( \|  u \|_{L^2(\Omega)}^2  + \|u_{yy} \|_{L^2(\Omega)}^2\Big)+L\Big(\|  u_t \|_{L^2(\Omega)}^2 + \|u^2u_x\|_{L^2(\Omega)}^2\Big) \nonumber \\
%===========================
+ \frac{2}{L^2} \int_{-B}^B u_x^2(0,y,t) \, dy \nonumber \\
%===========================
\leq (1+L)K_1   \Big(\frac{8}{L} +L +\frac{4(1+L)^2}{L^2} +\frac{K_2}{(1+L)K_1}\Big) e^{-\gamma_2 t} +\|u^2u_x\|_{L^2(\Omega)}^2 \nonumber \\
%===========================
\end{eqnarray}
%-------------------------------
For the latter right-hand norm we write
\begin{eqnarray}\label{estimativa traço uxx L 3}
 \|u^2u_x\|_{L^2(\Omega)}^2 \leq \sup_{(x,y)\in\Omega}u^4(x,y,t)\|u_x\|^2(t) \nonumber \\
%===========================
\leq \Big(6C^{10}_{N10}\|u_0\|^2\|\nabla u\|^8(t) + 2\|\nabla u\|^2(t)+2\| u_{xy}\|^2(t) \Big)\|u_x\|^2(t) \nonumber \\
%===========================
\leq \Big(6C^{10}_{N10}\|u_0\|^2C_{\|u_0\|}^4(1+L)^4K_1^4e^{-4\gamma_2 t} + 2C_{\|u_0\|}(1+L)K_1e^{-\gamma_2 t}+K_2 e^{-\gamma_2 t}\Big)\|u_x\|^2(t) \nonumber \\
%===========================
\leq \Big(6C^{10}_{N10}\|u_0\|^2C_{\|u_0\|}^4(1+L)^4K_1^4 + 2C_{\|u_0\|}(1+L)K_1+K_2 \Big)e^{-\gamma_2 t}\|u_x\|^2(t) \nonumber \\
%===========================
\leq \Big(6C^{10}_{N10}\|u_0\|^2C_{\|u_0\|}^4(1+L)^4K_1^4 + 2C_{\|u_0\|}(1+L)K_1+K_2 \Big)2(1+L)^3K_1e^{-\gamma_2 t} e^{-\gamma_2 t}\nonumber \\
%===========================
= K_3 e^{-2\gamma_2 t} \nonumber \\
%===========================
\leq K_3 e^{-\gamma_2 t}
\end{eqnarray}
%-------------------------------
with
\begin{eqnarray}\label{estimativa traço uxx L 4}
K_3=\Big(6C^{10}_{N10}\|u_0\|^2C_{\|u_0\|}^4(1+L)^4K_1^4 + 2C_{\|u_0\|}(1+L)K_1+K_2 \Big)2(1+L)^3K_1.
\end{eqnarray}
Finally,
\begin{eqnarray}\label{estimativa traço uxx L 5}
 \int_{-B}^B  u^2_{xx}(L,y,t)\, dy \leq   K_4 e^{-\gamma_2 t}
\end{eqnarray}
where
\begin{eqnarray}\label{estimativa traço uxx L 6}
K_4=(1+L)K_1   \Big(\frac{8}{L} +L +\frac{4(1+L)^2}{L^2} +\frac{K_2}{(1+L)K_1}\Big) + K_3.
\end{eqnarray}

\subsection{Estimate VII} Differentiate the equation with respect to

$x$, multiply by $(1+x) u_{x}$ and integrate over $\Omega.$ The
result is

\begin{eqnarray}\label{decaimento uxx}
%-------------------------------
\frac{d}{dt}\left( 1+x, u_x^2 \right) (t)+\|\nabla u_x \|^2(t) + 2\| u_{xx} \|^2(t)+ \int_{-B}^B \left[u_{xx}^2(0,y,t)+u^2_{xy} (0,y,t)\right]dy \nonumber\\
%-------------------------------
 =2\int_{-B}^B  \left[u_xu_{xxx}(0,y,t) -   u_xu_{xx}(0,y,t)+\frac{1}{2}u_{x}^2(0,y,t)+\frac{(1+L)}{2}u_{xx}^2(L,y,t)\right]dy \nonumber \\
 +\|u_x\|^2(t)-\frac{2}{3}\int_{\Omega} (1+x)u_x \big( u^3\big)_{xx}   \, d\Omega. \nonumber\\
 %-------------------------------
\end{eqnarray}
%-------------------------------
%-------------------------------
Inserting this into the equations gives

\begin{eqnarray}\label{decaimento uxx 1}
%-------------------------------
2\int_{-B}^B  u_xu_{xxx}(0,y,t)dy &=& -2\int_{-B}^B  u_x[u_t+u_x+u_{xyy}+u^2u_x](0,y,t)dy \nonumber\\
%-------------------------------
&=&-2\int_{-B}^B   u_x^2(0,y,t)dy-2\int_{-B}^B   u_x u_{xyy}(0,y,t)dy \nonumber\\
%-------------------------------
&=& -2\int_{-B}^B   u_x^2(0,y,t)dy + 2\int_{-B}^B  u_{xy}^2(0,y,t)dy .\nonumber\\
%-------------------------------
\end{eqnarray}
%-------------------------------
Substituting (\ref{decaimento uxx 1}) into (\ref{decaimento uxx})
provides

\begin{eqnarray}\label{decaimento uxx 2}
%-------------------------------
 \frac{d}{dt}\left( 1+x, u_x^2 \right) (t)+\|\nabla u_x \|^2(t) + 2\| u_{xx} \|^2(t)+ \int_{-B}^B \left[u_{xx}^2(0,y,t)+u^2_{x} (0,y,t)\right]dy \nonumber\\
%-------------------------------
 =2\int_{-B}^B  \left[ \frac{1}{2}u_{xy}^2(0,y,t)+\frac{(1+L)}{2}u_{xx}^2(L,y,t) -   u_xu_{xx}(0,y,t)\right]dy \nonumber \\
+\|u_x\|^2(t)-\frac{2}{3}\int_{\Omega} (1+x)u_x \big( u^3\big)_{xx}   \, d\Omega. \nonumber\\
 %-------------------------------
=2\int_{-B}^B  \left[ \frac{1}{2}u_{xy}^2(0,y,t)+\frac{(1+L)}{2}u_{xx}^2(L,y,t) -   u_xu_{xx}(0,y,t)\right]dy \nonumber \\
 +\|u_x\|^2(t)+2\int_{\Omega} \big[ u^2u_{x}^2+(1+x)u^2u_{x}u_{xx}\big]d\Omega.
\end{eqnarray}
%-------------------------------
%-------------------------------
%-------------------------------
As a consequence, we have
\begin{eqnarray}\label{decaimento uxx 3}
%-------------------------------
 \frac{d}{dt}\left( 1+x, u_x^2 \right) (0)+\|\nabla u_x \|^2(t) + \| u_{xx} \|^2(t)+ \frac{1}{2}\int_{-B}^B u_{xx}^2(0,y,t)dy \nonumber\\
%-------------------------------
 \leq \int_{-B}^B  \left[3u^2_{x} (0,y,t)+ u_{xy}^2(0,y,t)+(1+L)u_{xx}^2(L,y,t) \right]dy \nonumber \\
 +\|u_x\|^2(t) + (1+L)^2\|u^2 u_x\|^2(t)+2 \| u u_x \|^2(t).
\end{eqnarray}
%-------------------------------
%-------------------------------
%-------------------------------
Thanks to (\ref{decaimento ut 23}), (\ref{decaimento 6}),
(\ref{decaimento 6.1}), (\ref{estimativa traço uxx L 3} ) and
(\ref{estimativa traço uxx L 5} ), one concludes
\begin{eqnarray}\label{decaimento uxx 4}
%-------------------------------
\frac{d}{dt}\left( 1+x, u_x^2 \right) (t)+ \|\nabla u_x \|^2(t) + \| u_{xx} \|^2(t)+ \frac{1}{2}\int_{-B}^B u_{xx}^2(0,y,t)dy \nonumber\\
%-------------------------------
 \leq   \left[6(1+L)^3K_1+ K_2 +(1+L)K_4+2(1+L)^2K_1+(1+L)^2K_3 \right]e^{-\gamma_2 t} \nonumber \\
 + 2 \| u u_x \|^2(t)\nonumber\\
%-------------------------------
\leq   \left[6(1+L)^3K_1+ K_2 +(1+L)K_4+2(1+L)^2K_1+(1+L)^2K_3 \right]e^{-\gamma_2 t} \nonumber \\
  + 2\sup_{(x,y)\in\Omega}u^2(x,y,t) \|  u_x \|^2(t)\nonumber\\
%-------------------------------
\leq \cdots  + 4(1+L)^2K_1 \big(\|u\|^2(t) + \|\nabla u \|^2(t) +\|u_{xy}\|^2(t)  \big)e^{-\gamma_2 t}
\nonumber\\
%-------------------------------
\leq \cdots  + 4(1+L)^2K_1\big((1+L)\|u_0\|^2 + C_{\|u_0\|}(1+L)K_1 +K_2  \big)e^{-2\gamma_2 t}\nonumber\\
%-------------------------------
\leq K_5 e^{-\gamma_2 t}
\end{eqnarray}
%-------------------------------
where
\begin{eqnarray}\label{decaimento uxx 5}
%-------------------------------
 K_5 = 6(1+L)^3K_1+ K_2 +(1+L)K_4+2(1+L)^2K_1+(1+L)^2K_3 \\ + 4(1+L)^2K_1\big((1+L)\|u_0\|^2 + C_{\|u_0\|}(1+L)K_1 +K_2  \big).
\end{eqnarray}
%-------------------------------
Integrate (\ref{decaimento uxx 4}) in $t\in[0,T].$ The result reads
\begin{eqnarray}\label{decaimento uxx 5}
\| u_{xx} \|^2_{L^2_tL^2_{xy}}+
\left( 1+x, u_x^2 \right) (T)+ \|\nabla u_x \|^2_{L^2_tL^2_{xy}}
\frac{1}{2}\int_0^T \int_{-B}^B u_{xx}^2(0,y,t)\,dydt \nonumber\\
%-------------------------------
\leq K_5 \int_0^T e^{-\gamma_2 t} \,dt + \left( 1+x, u_x^2 \right) (0)\nonumber\\
%-------------------------------
\leq \frac{K_5}{\gamma_2}\Big(1-e^{-\gamma_2 T}  \Big) +
(1+L)\|u_{0x}\|^2.
\end{eqnarray}
%-------------------------------
Note that all the constants $K_i$ are proportional to
$\|u_0\|_{D(A)}.$ Since all the estimates do not depend upon the
$T_0,$ the local solution $u\in X_{T_0}$ can be continued for all
$T>0$ with the decay rate described above.

The uniqueness of solution $u\in X_T$ is proven by the usual way,
using similar computations as in lemma \ref{contraÃ§Ã£o}.

The proof of Theorem \ref{teorema principal} is completed.

%\section{Acknowledgments}
%We appreciate very much fruitful and constructive comments and
%suggestions of both of the Reviewers.

% N. Hongsit, M.A. Allen, G. Rowlands,
% Growth rate of transverse instabilities of solitary pulse solutions
% to a family of modified Zakharov-Kuznetsov equations,
% Physics Letters A, 372(14), 2420 (2008).

\end{document}